\documentclass{amsart}

\newcommand{\thmref}[1]{Theorem~\ref{#1}}
\newcommand{\secref}[1]{section~\ref{#1}}
\newcommand{\lemref}[1]{Lemma~\ref{#1}}
\newcommand{\propref}[1]{Prop\-o\-si\-tion~\ref{#1}}

\newcommand{\corref}[1]{Cor\-ol\-lary~\ref{#1}}

\newcommand{\tr}[1]{\hspace{-1ex}\phantom{x}^t #1}
\newcommand{\defnref}[1]{Def\-i\-ni\-tion~\ref{#1}}

\DeclareMathOperator{\Stab}{Stab}
\DeclareMathOperator{\Sym}{Sym}

\renewcommand{\a}{\mathbf{a}}
\renewcommand{\b}{\mathbf{b}}

\renewcommand{\r}{\mathbf{r}}

\newcommand{\cS}{\mathcal{S}}
\newcommand{\cSkn}{\cS_k(\Gamma^n)}

\newcommand{\diag}{\text{diag}}

\newcommand{\bbZ}{\mathbb{Z}}
\newcommand{\bbR}{\mathbb{R}}
\newcommand{\bbQ}{\mathbb{Q}}

\newcommand{\QW}{\bbQ[x_0,\dots,x_n]^{W_n}}

\newcommand{\C}{\mathbb{C}}

\newcommand{\bleq}{\boldsymbol{\leq}}

\newcommand{\cHint}{\underline{\mathcal{H}}_{n,p}}
\newcommand{\cHp}{\mathcal{H}_{n,p}}

\newcommand{\CWn}{(\bbC^\times)^{n+1}/W_n}
\newcommand{\bbC}{\mathbb{C}}
\newcommand{\Hom}{\text{Hom}}

\newcommand{\Orbit}{\text{Orbit}}
\newcommand{\cH}{\mathcal{H}}

\newtheorem{theorem}{Theorem}[section]
\newtheorem{thm}[theorem]{Theorem}

\newtheorem{lem}[theorem]{Lemma}

\newtheorem{prop}[theorem]{Proposition}
\newtheorem{cor}[theorem]{Corollary}

\theoremstyle{definition}
\newtheorem{definition}[theorem]{Definition}
\newtheorem{defn}[theorem]{Definition}

\newtheorem{notation}[theorem]{Notation}
\theoremstyle{remark}

\numberwithin{equation}{section}

\begin{document}

\title{Computing the Satake $p$-parameters of Siegel modular forms}

\author{Nathan C. Ryan}
\address{Department of Mathematics, Dartmouth College, Hanover, NH 03755}
\email{nathan.c.ryan@dartmouth.edu}
\thanks{The first author was supported in part by a GAANN Fellowship.}



\date{October 25, 2004 and, in revised form, \today.}



\begin{abstract}
Hecke eigenvalues of
  classical modular forms often encode a wealth of arithmetic data.
  The Satake $p$-parameters of a Siegel modular form play a role
  analogous to the one played by Hecke eigenvalues in the
  characterization of classical modular forms.  In this paper we
  present an algorithm by which we can compute the Satake
  $p$-parameters of a Siegel modular form $F$ if we are given the
  Hecke eigenvalues of $F$ with respect to the generators of the local
  Hecke algebra.  We compute explicit examples and relate our
  computations to the Ramanujan-Petersson Conjecture.
\end{abstract}

\maketitle



In the 1960's, Satake \cite{Satake} showed the
existence of an isomorphism $\Omega$ between the local, integral Hecke
algebra $\cHint$ and $\QW$ (see Notation and Terminology for
precise definitions), a
ring of polynomials with rational coefficients that are invariant under
the action of a signed permutation group $W_n$  ($W_n$ is the Weyl
group associated to the symplectic group over a local field).  Proving the existence of
$\Omega$ is a key step in proving the existence of the Satake isomorphism
$\Psi:\text{Hom}_{\bbC}(\cHp,\bbC)\to (\bbC^\times)^{n+1}/W_n$ \cite{Abook}.

Let $F$ be a Hecke eigenform of degree $n$ and weight $k$ and let the map
$\lambda\in\text{Hom}_{\bbC}(\cHp,\bbC)$ be given by the rule $T\mapsto
\lambda_F(T)$ (here $\lambda_F(T)$ is the eigenvalue of $F$ with
respect to the operator $T$).  Then the entries of the $(n+1)$-tuple
$(\alpha_{0,p},\dots,\alpha_{n,p}) (\text{mod } W_n)$ of complex numbers $\Psi(\lambda)$
are called the Satake $p$-parameters of $F$.

Satake $p$-parameters are to Siegel modular forms
and their associated $L$-functions what Hecke eigenvalues are to
classical modular forms and their associated $L$-functions.  Satake
$p$-parameters are used to define $L$-functions, to generalize the
Ramanujan-Petersson Conjecture, and to develop results
suggestive of multiplicity-one theorems \cite{Veenstra}.  Because of
their prevalence, it is certainly of interest to be able to compute
the Satake parameters of a given modular form $F$.  

From the Satake
$p$-parameters of a Hecke eigenforms $F$, we can determine the Hecke
eigenvalue of $F$ with respect to any Hecke operator $T$ by means of
the spherical map $\Omega$ \cite{Satake},\cite{Gross}.  The spherical map \cite{Abook} is defined on right cosets
whereas Hecke operators are linear combinations of double cosets.  The
image of a double coset in the Hecke algebra is the sum of the image
of the right cosets into which the double coset decomposes.  Two questions naturally arise:
(1)  Can
we explicitly determine a matrix representation of the spherical map?
and (2) What Hecke data for $F$ do we need to determine the
Satake $p$-parameters of $F$?

Our starting point is the following proposition of Andrianov \cite{Abook}, 
\begin{prop}\label{Aprop}
Let $p$ be prime.  Each non-zero
$\lambda\in\text{Hom}_{\bbC}(\cHp,\bbC)$ has the following form:  for $T\in \cHp$,
\[
T\mapsto
\lambda_A(T)=\Omega(T)|_{(x_0,\dots,x_n)=(\alpha_{0,p},\dots,\alpha_{n,p})}
\]
where $A=(\alpha_{0,p},\dots,\alpha_{n,p})$ is a vector of non-zero complex
numbers depending only on $\lambda$.
\end{prop}
\noindent Thus, to compute the image of $\Psi$, we first need to compute the
image of $\Omega$.  A paper by Krieg \cite{Krieg} tells us how to
compute the $\Omega$-image of the generators of $\cHint$ which is all
that is needed for our work.

An outline of this paper is as follows:  using Proposition \ref{Aprop}
we get a system of $n+1$ polynomial equations in $n+1$ unknowns.
Performing some straightforward algebraic manipulation, we change this
into a system of $n$ polynomial equations in $n$ unknowns.  The
algorithm we develop takes as input this manipulated system and
outputs a representative of the Satake $p$-parameter $W_n$-equivalence
class.  Using this algorithm we verify the (strong) Ramanujan
Petersson Conjecture for particular Hecke eigenforms based on Hecke
eigenvalue computations by Skoruppa \cite{Skoruppa} and Breulmann and
Kuss \cite{Breul2}.

\specialsection*{NOTATION AND TERMINOLOGY}

Our notation follows closely that of Andrianov \cite{Abook}.  
\begin{itemize}
\item by $\Gamma^n$ we denote Siegel's modular group $Sp_n(\bbZ)$ and
  by $G$ we denote $GSp^+_n(\bbR)=\left\{\alpha=\left(\begin{smallmatrix}A&B\\C&D\end{smallmatrix}\right)\in
  M_{2n}(\bbR): \alpha\text{ satisfies (i) and (ii)}\right\}$ where (i) is $\tr A D-\tr CB=r(\alpha)I_n$
for some 
  $r(\alpha)\in\bbR^+$ and (ii) is that $\tr AC,\tr BD$ are
  symmetric.  The positive real number $r(\alpha)$ is the
  \textit{similitude of $\alpha$}.
\item  let $R(\Gamma,S)$ be the space of finite formal $\bbZ$-linear
  combinations of right cosets of the form $\Gamma\alpha$ where
  $\alpha\in S=GSp^+_n(\bbQ)$.  By $R(\Gamma,S)^\Gamma$ we denote the
  stabilizer of $R(\Gamma,S)$ under the action given by right
  multiplication by elements of $\Gamma$.  Then
  $\cH_n=\otimes_\bbZ\bbC$ denotes the \textit{global Hecke algebra}.
  For $S_p=S\cap M_{2n}(\bbZ[p^{-1}])$ we define
  $R(\Gamma,S_p)^\Gamma$ in a similar way.  Then
  $\cHp=R(\Gamma,S_p)^\Gamma\otimes_\bbZ\bbC$ denotes the \textit{local Hecke
    algebra}.  Finally, for $\underline{S}_p=S_p\cap M_{2n}(\bbZ)$ we
  define $R(\Gamma,\underline{S}_p)^\Gamma$ in a similar way and
  denote by $\cHint=R(\Gamma,\underline{S}_p)^\Gamma\otimes_\bbZ\bbC$
  the \textit{integral Hecke algebra}.

\item by $W_n$ we denote the group of $\bbQ$-automorphisms of
  $\bbQ(x_0,x_1^{\pm 1},\dots,x_n^{\pm 1})$ generated by the
  permutations of the variables $x_1,\dots,x_n$ and for $1\leq i\leq
  n$ the maps $\tau_i$ given by
\[
\tau_i(x_0)=x_0x_i,\,\tau_i(x_i)=x_i^{-1},\text{ and
}\tau_i(x_j)=x_j\,\,(j\neq i,0).
\]
We set $\bbQ[x_1,\dots,x_n]=\left\{f\in\bbQ[x_1,\dots,x_n]:\sigma(f)=f\text{ for all
  }\sigma\in S_n\right\}$ and $\QW=\left\{f\in \bbQ[x_0,\dots,x_n]:\sigma(f)=f\text{ for all
  }\sigma\in W_n\right\}$.  Following Andrianov, we call the map
$\Psi:\Hom(\cHp,\bbC)\to(\bbC^\times)^{n+1}/W_n$ the \textit{Satake
  isomorphism} and the map $\Omega:\cHint\to\QW$ the
\textit{spherical map}.
\item by $S_n$ we denote the symmetric group on $n$ letters and it
  acts on $f\in\bbQ[x_1,\dots,x_n]$ by permuting the subscripts of the
  indeterminates $x_i$.  Then for $f\in\bbQ[x_1,\dots,x_n]$ we define the \textit{$S_n$ symmetrized
    polynomial of $f$} to be
\[
\Sym_{S_n}(f)=\sum_{\sigma\in S_n/\Stab(f)}\sigma(f)
\]
and for $f\in\QW$ we define the \textit{$W_n$ symmetrized polynomial of $f$} to be
\[
\Sym_{W_n}(f)=\sum_{\sigma\in W_n/\Stab(f)}\sigma(f).
\]
For brevity we will denote $\Sym_{W_n}(f)$ by $[f]$.
\item Let $\b=(b_1,\dots,b_n)\in\bbZ^n$.  Then by $p^\b$ we mean the
  matrix $\diag(p^{b_1},\dots,p^{b_n})$ and by $x^\b$ we mean the
  monomial $x_1^{b_1}\cdots x_n^{b_n}$.
\item by $\cSkn$ we denote the space of Siegel cusp forms of weight
  $k$ and degree $n$; by a \textit{Hecke eigenform} we mean a cusp
  form that is an eigenforms for all the Hecke operators $T(m)$.
\item there is an imbedding $\Gamma\alpha\Gamma=\cup_{i=1}^\nu
  \Gamma\alpha_i\mapsto \sum_{i=1}^\nu\Gamma\alpha_i$ of double cosets
  into $\cH_n$.  It makes sense, then, to say $\cHint$ is generated
  by the $n+1$ double cosets $T_0(p)=\Gamma\diag(I_n;pI_n)\Gamma$ and
  for $1\leq i \leq n$
  $T_i(p^2)=\Gamma\diag(I_i,pI_{n-i};p^2I_i,pI_{n-i})\Gamma$.  We also
  define the extra operator $T_0(p^2)=(I_n;p^2I_n)$.  
\item $\Gamma\alpha\Gamma=\cup_{i=1}^\nu$ acts on
  $\cSkn$ by
  $F|_k\Gamma\alpha\Gamma=r(\alpha)^{nk-n(n+1)}\sum_{i=1}^\nu
  F|_k(\alpha_i)$ where $r(\alpha)$ is the similitude of $\alpha$ and 
\[
F_k(\alpha_i)=\det(CZ+D)^{-k}F((AZ+B)(CZ+D)^{-1})\text{ where }
\alpha_i=\left ( \begin{smallmatrix} A&B\\C&D\end{smallmatrix}\right).
\]
By $\lambda_F(T)$ we mean the eigenvalue of $F$ with respect to
operator $T\in\cH$.
\item Let $F\in\cSkn$ be a simultaneous eigenform.  The \textit{Satake $p$-parameters associated
  to $F$} is the $(n+1)$-tuple $(\alpha_{0,p},\dots,\alpha_{n,p})\in
  \CWn$ which is the image of the homomorphism $T\mapsto \lambda_F(T)$
  under the Satake isomorphism $\Psi:\Hom_\bbC(\cHp,\bbC)\to
  \CWn$.
\end{itemize}

Before starting in earnest, we state and prove the following easy lemma:
\begin{lem}\label{orbits}
Let $\a=(a_1,\dots,a_n)\in\bbZ_+^n$ and $r\geq 0$.  Under the action of $W_n$ we obtain the
following orbit:
\begin{align*}
\Orbit_{W_n}(x_0^r x^\a)&=\Orbit_{S_n}(\{x_0^rx_1^{\epsilon_1}\cdots
x_n^{\epsilon_n}:\epsilon_i=a_i \text{ or } r-a_i\}\\
&=\left\{x_0^r x_{\sigma(1)}^{\epsilon_{\sigma(1)}}\cdots
  x_{\sigma(n)}^{\epsilon_{\sigma(n)}}:\sigma\in S_n,\epsilon_i=a_i \text{ or } r-a_i\right\}.
\end{align*}
\end{lem}
\begin{proof}  We see the first equality directly:  applying $\tau_i$
  to $x_0^r x^\a$ we get
\begin{align*}
\tau_i(x_0^rx^\a)&=x_0^rx_i^rx_1^{a_1}\cdots x_i^{-a_i} \cdots x_n^{a_n}\\
&=x_0^rx_1^{a_1}\cdots x_i^{r-a_i} \cdots x_n^{a_n}.
\end{align*}
Since $W_n$ is generated by the $\tau_i$ and $\sigma\in S_n$ we get
the first equality.  The second equality follows immediately.
\end{proof}

\section{A Result of Krieg}\label{Kriegsec}  Let $\r=(2,\dots,2)\in\bbZ^n$.  Write the operator $T_i(p^2)$ ($0\leq
i\leq n$) as $\Gamma\diag(p^{\r-\b_i};p^{\b_i})\Gamma$ so that
$\b_i=(\underbrace{2,\dots,2}_{n-i},\underbrace{1,\dots,1}_i)$.  In
\cite{Ryan1} we show that 
\begin{equation}\label{sys0}
\Omega(T_i(p^2))=\sum_{j=i}^n c(i,j)\Sym_{W_n}(x_0^2x^{\b_j})
\end{equation}
\noindent where a method to compute the coefficients is described.  We
show, in particular, that $c(j,j)$ is non zero.
Krieg (Corollary 2 in \cite{Krieg}) does one better than this and gives the following
explicit expressions for the coefficients $c(i,j)$ (notice that we
normalize our slash operator differently, so we divide his
coefficients by $p^{nk-n(n+1)/2}$ \cite{Abook}):
\begin{prop}\label{Krieg}
Let $i,j\geq 0$.  Then
\begin{align*}
c(2i+1,j)&=(p^{j+1}-1)c(2i,j+1)\\
c(2i,j)&=p^{-nk+n(n+1)/2}\sum_{l=0}^{i}
(-1)^l\binom{2i+j}{i-l}p^{4i^2+4ij+2i+\frac{1}{2}j^2+\frac{1}{2}j-2jl-l^2}
  d(l,j)
\end{align*}
\noindent where 
\begin{align*} 
d(0,j)&:=1\\
d(1,j)&:=p^j+\frac{p^{2j+2}-1}{p^2-1}
\end{align*}
and 
\[
d(l,j):=\prod_{t=1}^l
\frac{p^{2j+2t}-1}{p^{2t}-1}+p^j\prod_{t=1}^{l-1}
\frac{p^{2j+2t}-1}{p^{2t}-1}
\]
for $l\geq 2$.
\end{prop}

\section{Notation for the Algorithm}\label{alg-not}

Suppose $[\Omega_2]$ is the upper triangular matrix whose $(i,j)$th entry
($j\geq i$) is the coefficient $c(i,j)$ from \propref{Krieg}.    

\propref{Aprop} relates the eigenvalues for $T\in\cHint$ with Satake
parameters.  Let $\lambda_{F,0}(p)$ denote the eigenvalue of a
simultaneous eigencuspform $F$ with respect to the Hecke operator
$T_0(p)$ and let $\lambda_i$ denote the eigenvalue of $F$ with respect to
$T_i(p^2)$ ($0\leq i \leq n$).  According to equation \ref{sys0} and \propref{Aprop} we get the following system of
equations (see Notation above):
\begin{align}\label{sys1}
\lambda_{F,0}(p)&=x_0(1+x_1)\cdots(1+x_n)\nonumber\\
\lambda_n&=c_{11}[x_0^2x^{(1,\dots,1)}]\nonumber\\
\lambda_{n-1}&=c_{12}[x_0^2x^{(1,\dots,1)}]+c_{22}[x_0^2x^{(2,1,\dots,1)}]\\
\lambda_{n-2}&=c_{13}[x_0^2x^{(1,\dots,1)}]+c_{23}[x_0^2x^{(2,1,\dots,1)}]+c_{33}[x_0^2x^{(2,2,1,\dots,1)}]\nonumber\\
&\vdots\nonumber\\
\lambda_0&=c_{1,n+1}[x_0^2x^{(1,\dots,1)}]+
\cdots+c_{n+1,n+1}[x_0^2x^{(2,\dots,2)}].\nonumber
\end{align}
\noindent We have left the equations in polynomial notation rather than
in Satake parameter notation but we retain the assumption that the
$x_i$ cannot take 0 as a value. 

We note that the
system is ``upper triangular'' in the sense that all but the last term of the right hand
side of each equation is the same as the right hand side of the
preceding equation.  Hence this system becomes, for constants
$k_1,\dots k_{n+1},$ (note that we drop the first equation)
\begin{align}\label{sys2}
\lambda_n&=c_{11}[x_0^2x^{(1,\dots,1)}]:=c_{11}k_1\nonumber\\
\lambda_{n-1}&=c_{12}k_1+c_{22}[x_0^2x^{(2,1,\dots,1)}]:=c_{12}k_1+c_{22}k_2\nonumber\\
\lambda_{n-2}&=c_{13}k_1+c_{23}k_2+c_{33}[x_0^2x^{(2,2,1,\dots,1)}]:=c_{13}k_1+c_{23}k_2+c_{33}k_3\\
&\vdots\nonumber\\
\lambda_0&=\sum_{i=1}^{n} c_{i,n+1}k_i+
c_{n+1,n+1}[x_0^2x^{(2,\dots,2)}]:=\sum_{i=1}^{n} c_{i,n+1}k_i+
c_{n+1,n+1}k_{n+1}.\nonumber
\end{align}
 For the sake of clarity in the algorithm, we need to rearrange and rewrite the
 equations.  But first we need to introduce the following piece of
 notation and prove the following lemma:
\begin{notation}  Let $k,j\in\bbZ$ and $n\in\bbZ^+$.  Let $\mathfrak{S}^{(n)}_{k,j}=\Sym_{S_n}(x_1^2\cdots
 x_k^2 x_{k+1}\dots x_{k+j}$).  If $k+1<0$ or $k+j>n$,
 $\mathfrak{S}^{(n)}_{k,j}=0$.  
\end{notation}
\noindent In words, $\mathfrak{S}^{(n)}_{k,j}$ is the
 $S_n$-symmetrized polynomial of a monomial in $\bbQ[x_1,\dots,x_n]$
 with $k$ of the $x_i$ squared and $j$ of the $x_i$ to the first power.

The following lemma relates $S_n$ symmetrized polynomials and the
$W_n$ symmetrized polynomials in system \ref{sys2}:

\begin{lem}\label{frakS}  Let $n\in\bbZ^+$ and the integer  $0\leq j\leq n$.  Then
\begin{equation}\label{lem-eq}
[x_0^2x_1^2\dots x_{n-j}^2
x_{n-j+1}\dots x_n]=x_0^2\sum_{k=0}^{n-j}\mathfrak{S}^{(n)}_{k,j}.
\end{equation}
\end{lem}
\begin{proof}
As both sides of equation \ref{lem-eq} are polynomials where every
monomial has coefficient one, we argue that the set of terms of the
right hand side of equation \ref{lem-eq} is contained in the set of
terms on the left side of the equation and that these two sets have
the same size.  

Let $x_0^2x^\a$ be a monomial on the left hand side of the equation.
We know that $x_0^2 x^\a\in \Orbit_{W_n}(x_0^2x_1^2\cdots
x_{n-j}^2x_{n-j+1}\cdots x_n)$.  By
\lemref{orbits} we know that 
\begin{multline}\label{lem-eq2}
\Orbit_{W_n}(x_0^2x_1^2\cdots x_{n-j}^2x_{n-j+1}\cdots
x_n)=\\
\{x_0^2 x_{\sigma(1)}^{\epsilon_{\sigma(1)}}\cdots
  x_{\sigma(n)}^{\epsilon_{\sigma(n)}}:\sigma\in S_n,
    \epsilon_1=\cdots=\epsilon_j=1,\epsilon_{j+1},\dots,\epsilon_{n}=0
    \text{ or } 2\}
\end{multline}
Suppose, as the left hand side of equation \ref{lem-eq2} is invariant
under permutations of $x_1,\dots,x_n$ , that
$\epsilon_1,\dots,\epsilon_k$ are 2,
$\epsilon_{k+1},\dots,\epsilon_{k+j}$ are 1 and
$\epsilon_{k+j+1},\dots,\epsilon_n$ are 0. Then the monomial $x_0
x^\a$ is a term in $\Sym_{S_n}(x_1^2\dots x_k^2 x_{k+1}\cdots
x_{k+j})=\mathfrak{S}^{(n)}_{k,j}$.

We count the number of terms on each side of equation \ref{lem-eq}.
According to \lemref{orbits} there are $\binom{n}{j}$ places for
$\epsilon_i$ to equal 1 and the remaining $n-j$ $\epsilon_i$'s are
either 0 or 2.  So the total number of terms on the left hand side of
equation \ref{lem-eq} is $\binom{n}{j} 2^{n-j}$.

For each $0\leq k \leq n-j$ the polynomial $\mathfrak{S}^{(n)}_{k,j}$
has $n(n-1)\cdots(n-k+1)$ places to let $\epsilon_i$ be 2 and since
order does not matter we get that are $n(n-1)\cdots(n-k+1)/k!$ way to
arrange the 2's.  Since there are $n-k$ remaining $\epsilon_i$'s and
$j$ of these have to be 1's and the rest 0's, we have a total of
$\binom{n-k}{j} n(n-1)\cdots(n-k+1)/k!$ terms in
$\mathfrak{S}^{(n)}_{k,j})$.  A straightforward computation shows:
\begin{align*} 
\sum_{k=0}^{n-j} \binom{n-k}{j} 
\frac{n(n-1)\cdots(n-k+1)}{k!}
= \binom{n}{j} 2^{n-j}.
\end{align*}\end{proof}

We apply \lemref{frakS} to change the
equations in system \ref{sys2} to sums of $S_n$-invariant polynomials:
\begin{align}\label{sys3}
k_1=\lambda_{n}/c_{11}&=x_0^2\mathfrak{S}^{(n)}_{0,n}x_1\cdots
x_n=x_0^2x_1\cdots x_n\nonumber\\
k_2=\frac{\lambda_{n-1}-c_{12}k_1}{c_{22}}&=x_0^2\sum_{k=0}^1\mathfrak{S}^{(n)}_{k,n-1}
\nonumber\\
k_3=\frac{\lambda_{n-2}-c_{13}k_1-c_{23}k_2}{c_{33}}&=x_0^2\sum_{k=0}^2\mathfrak{S}^{(n)}_{k,n-2}\nonumber\\
&\vdots\nonumber\\
k_j=\frac{\lambda_{n-j}-\sum_{i=1}^{j-1} c_{ij}k_i)}{c_{jj}}&=x_0^2\sum_{k=0}^{j-1}\mathfrak{S}^{(n)}_{k,n-j+1}\\
&\vdots\nonumber\\
k_{n+1}=\frac{\lambda_{0}-\sum_{i=1}^{n} c_{i,n+1}k_i)}{c_{n+1,n+1}}&=x_0^2\sum_{k=0}^n\mathfrak{S}^{(n)}_{k,0}\nonumber
\end{align}

The first equation in system \ref{sys3} says that the monomial
$x_0^2x_1\cdots x_n=k_1$ is constant.  We eliminate the variable $x_0$ in
system \ref{sys3} by multiplying both sides by $x_1\cdots x_n$ and
subsequently divide by $k_1$ to get
the following system of $n$ equations in $n$ unknowns:
\begin{align}\label{sys3a}
\frac{k_2}{k_1}x_1\cdots x_n&=\sum_{k=0}^1\mathfrak{S}^{(n)}_{k,n-1}\nonumber\\
\frac{k_3}{k_1}x_1\cdots x_n&=\sum_{k=0}^2\mathfrak{S}^{(n)}_{k,n-2}\nonumber\\
&\vdots\nonumber\\
\frac{k_j}{k_1}x_1\cdots x_n&=\sum_{k=0}^{j-1}\mathfrak{S}^{(n)}_{k,n-j+1}\\
&\vdots\nonumber\\
\frac{k_{n+1}}{k_1}x_1\cdots x_n&=\sum_{k=0}^n\mathfrak{S}^{(n)}_{k,0}\nonumber
\end{align}
By defining new constants $c_j$ to be $\frac{k_{n+2-j}}{k_1}$ and
doing some rearranging we get
the system the algorithm takes as input:
\begin{align}\label{sys4}
f_1&:=\sum_{k=0}^{n} \mathfrak{S}^{(n)}_{k,0}-c_1x_1\dots x_n =0\nonumber\\
f_2&:=\sum_{k=0}^{n-1} \mathfrak{S}^{(n)}_{k,1} -c_2x_1\dots x_n=0\nonumber\\
&\vdots\nonumber\\
f_i&:=\sum_{k=0}^{n-i+1}\mathfrak{S}^{(n)}_{k,i-1}- c_ix_1\dots x_n=0\\
&\vdots\nonumber\\
f_n&:= \sum_{k=0}^{1}\mathfrak{S}^{(n)}_{k,n-1}-c_n x_1\dots x_n=0\nonumber.
\end{align}

To solve system \ref{sys4}, one more piece of notation is necessary
\begin{defn}\label{Spoly}  Let $f,g\in\bbQ[x_1,\dots,x_n]$ so that the
  leading term of $f$ is $x^\a=x_1^{a_1}\cdots x_n^{a_n}$ and
  the leading term of $g$ is $x^\b=x_1^{b_1}\cdots
  x_n^{b_n}$.  If $\a\bleq \b$ then we define the \textit{$S$-polynomial of
  $f$ and $g$} to be 
\[
S(f,g)=\frac{x^\b}{x^\a} f - g.
\]
\end{defn}
\noindent The main point of \defnref{Spoly} is to have a way to cancel off the leading
terms, with respect to the lexicographic order, of polynomials $f$ and $g$.  Notice that adding $S(f,g)$ to
a system of polynomials that contains $f$ and $g$ does not affect the
solution set to the system of polynomials.  See, for examples,
\cite{CoxLittleOshea}, for this and other techniques from algebraic geometry.

\section{The Algorithm}
With the notation from the previous section, we are ready to
present the algorithm to compute Satake parameters given the Hecke
eigenvalues for $T_0(p)$ and $T_i(p^2)$ for $0\leq i\leq n$.  Before
we present our algorithm for finding the solutions to
the system of polynomial equations $f_1=0,\dots,f_n=0$ as defined in
system \ref{sys4}, we make the following definition:

\begin{defn}  Given a polynomial $f\in\bbC[x_1,\dots,x_n]$ of degree
  $n$, the \textit{reciprocal $f^r(x)$ of $f(x)$} is $f^r(x)=x^nf(1/x)$.  If
  $f(x)=f^r(x)$, we call it a \textit{palindromic polynomial}.
\end{defn}

\begin{thm}\label{algorithm} The following algorithm $\mathcal{A}$ solves the systems of polynomial
  equations $f_1=0,\dots,f_n=0$ as defined in system \ref{sys4}.
\begin{enumerate} \item $f_{n+1}=S(f_1,f_2)$ and
  $f_{n+2}:=S(f_2,f_{n+1})$.
\item  $i:=3$; for $j:=3$ to $n$ do
\begin{list}{\setlength{\leftmargin}{.5in}
             \setlength{\rightmargin}{.5in}}
\item for $k:=1$ to $j$ do
\item $f_{n+i}:=S(f_{n+i-1},f_j)$;
\item $i:=i+1$
\item od;
\item $G^{(n)}_j:=f_{n+i-1}$
\end{list}
od;
\end{enumerate}
Algorithm $\mathcal{A}$ results in $G^{(n)}_n$ having degree
$(1,1,\dots,1,2n)$.  Moreover each term of $G^{(n)}_n$ is
divisible by $x_1\cdots x_{n-1}$ and the polynomial
$G^{(n)}_n/x_1\cdots x_{n-1}$ is palindromic.   
\end{thm}

We prove \thmref{algorithm} by the following lemmas and leave sample
applications of the algorithm $\mathcal{A}$ to the next section (to
see the algorithm in action, we suggest you look there first).  First we make the following definition:
\begin{defn} A finite sequence of numbers $a(k)$ for $k$ from $i$ to
  $j$ ($i,j,k\in\bbZ$) is palindromic if $a(k)=a(j-k+i)$.
\end{defn}

\begin{lem}\label{lem1} Suppose $i\in\bbZ$ and $1\leq i\leq n$.  Let $f_i$ be as
  in the statement of \thmref{algorithm}.  Then
\begin{multline*}
f_i:=(x_n^2+1)\left( (x_{n-1}^2+1)\sum_{k=0}^{n-i-1}
  \mathfrak{S}^{(n-2)}_{k,i-1}+x_{n-1}\sum_{k=0}^{n-i}\mathfrak{S}^{(n-1)}_{k,i-2}\right)\\
+x_n\left((x_{n-1}^2+1)\sum_{k=0}^{n-i}\mathfrak{S}^{(n-2)}_{k,i-2}+x_{n-1}(\sum_{k=0}^{n-i+1}\mathfrak{S}^{(n-2)}_{k,i-3}-c_ix_1\cdots
  x_{n-2})\right).
\end{multline*}
\end{lem}
\begin{proof}  By the definition of $f_i$ in system \ref{sys4} we have 
\[
f_i:=\sum_{k=0}^{n-i+1}\mathfrak{S}^{(n)}_{k,i-1}- c_ix_1\dots x_n.
\]
We parse $\mathfrak{S}^{(n)}_{k,i-1}$
according to the power of $x_n$ that appears in the monomials:
\[
\mathfrak{S}^{(n)}_{k,i-1}=x_n^2
\mathfrak{S}^{(n-1)}_{k-1,i-1}+x_n\mathfrak{S}^{(n-1)}_{k,i-2}+\mathfrak{S}^{(n-1)}_{k,i-1}
\]
\noindent The first summand has $x_n^2$ factored out so what remains
has only $k-1$ variables squared.  The second summand has $x_n$
factored out, so there are $i-2$ variables to the first.  All three
summands are superscripted by $n-1$ since we have accounted for the
appearance of $x_n$.  

Now
\begin{align*}
\sum_{k=0}^{n-i+1}\mathfrak{S}^{(n)}_{k,i-1}&- c_ix_1\dots
x_n\\
&=x_n^2\sum_{k=0}^{n-i+1}\mathfrak{S}^{(n-1)}_{k-1,i-1}+x_n\sum_{k=0}^{n-i+1}\mathfrak{S}^{(n-1)}_{k,i-2}+\sum_{k=0}^{n-i+1}\mathfrak{S}^{(n-1)}_{k,i-1}-c_ix_1\cdots
x_n\\
&=x_n^2\sum_{k=1}^{n-i+1}\mathfrak{S}^{(n-1)}_{k-1,i-1}+x_n
\sum_{k=0}^{n-i+1}\mathfrak{S}^{(n-1)}_{k,i-2}
+\sum_{k=0}^{n-i}\mathfrak{S}^{(n-1)}_{k,i-1}-c_ix_1\cdots x_n
\end{align*}
\noindent where the indexing of the first and second summand change
since $\mathfrak{S}^{(n-1)}_{-1,i-1}$ and
$\mathfrak{S}^{(n-1)}_{n-i+1,i-1}$ are zero.  Then, shifting the index
in the first summand we get 
\begin{equation}\label{f-i-1}
f_i:=(x_n^2+1)\sum_{k=0}^{n-i}
\mathfrak{S}^{(n-1)}_{k,i-1}+x_n(\sum_{k=0}^{n-i+1}\mathfrak{S}^{(n-1)}_{k,i-2}-c_ix_1\cdots
x_{n-1}).
\end{equation}
Next we apply a completely analogous argument to the
$\mathfrak{S}^{(n-1)}_{k,i}$ in equation \ref{f-i-1}, but we parse the
summands of
$\mathfrak{S}^{(n-1)}_{k,i}$ according to the power of $x_{n-1}$.
Hence, 
\begin{multline}\label{f-i-2}
f_i:=(x_n^2+1)\left( (x_{n-1}^2+1)\sum_{k=0}^{n-i-1}
  \mathfrak{S}^{(n-2)}_{k,i-1}+x_{n-1}\sum_{k=0}^{n-i}\mathfrak{S}^{(n-1)}_{k,i-2}\right)\\
+x_n\left((x_{n-1}^2+1)\sum_{k=0}^{n-i}\mathfrak{S}^{(n-2)}_{k,i-2}+x_{n-1}(\sum_{k=0}^{n-i+1}\mathfrak{S}^{(n-2)}_{k,i-3}-c_ix_1\cdots
  x_{n-2})\right).
\end{multline}\end{proof}

In equation \ref{f-i-2} we parse the monomials of $f_i$ according to
the power of $x_n$ and $x_{n-1}$ that appear.  Since $f_i$ is a
symmetric polynomial, if we apply a permutation $\sigma\in S_n$ so
that $n$ and $n-1$ get sent to $1$ and $n$ respectively, a priori,
equation \ref{f-i-2} would become:
\begin{multline}\label{f-i-3}
f_i:=(x_1^2+1)\left( (x_{n}^2+1)\sum_{k=0}^{n-i-1}
  \mathfrak{S}^{(n-2)}_{k,i-1}+x_n\sum_{k=0}^{n-i}\mathfrak{S}^{(n-1)}_{k,i-2}\right)+\\
x_n\left((x_{n}^2+1)\sum_{k=0}^{n-i}\mathfrak{S}^{(n-2)}_{k,i-2}+x_{n}(\sum_{k=0}^{n-i+1}\mathfrak{S}^{(n-2)}_{k,i-3}-c_ix_2\cdots
  x_{n-1})\right).
\end{multline}
\noindent As written, equation \ref{f-i-3} does not make sense since
  $\mathfrak{S}^{(n-2)}_{k,j}$ is a polynomial in $x_1,\dots,x_{n-2}$
  and thus there are too many $x_1$'s in the terms of $f_i$.  Thus we
  introduce the notation $\underline{\mathfrak{S}}^{(n-2)}_{k,j}$ to
  be the polynomial
  $\Sym_{S_{n-2}}(x_2^2,\dots,x_{k+1}^2,x_{k+2},\dots,x_{k+j+1})$ where
  the action of $S_{n-2}$ is on the indeterminates
  $x_2,\dots,x_{n-1}$.  With this notation, the correct version of
  equation \ref{f-i-3} is
\begin{multline}\label{f-i-4}
f_i:=(x_1^2+1)\left( (x_{n}^2+1)\sum_{k=0}^{n-i-1}
  \underline{\mathfrak{S}}^{(n-2)}_{k,i-1}+x_n\sum_{k=0}^{n-i}\underline{\mathfrak{S}}^{(n-2)}_{k,i-2}\right)+\\
x_n\left((x_{n}^2+1)\sum_{k=0}^{n-i}\underline{\mathfrak{S}}^{(n-2)}_{k,i-2}+x_{n}(\sum_{k=0}^{n-i+1}\underline{\mathfrak{S}}^{(n-2)}_{k,i-3}-c_ix_2\cdots
  x_{n-1})\right).
\end{multline}

\begin{lem}\label{lem2}  For $1\leq i\leq n$, let $G^{(n)}_i$ be as in the
  statement of \thmref{algorithm}.  If 
\begin{multline*}
G^{(n)}_i=x_1^2\sum_{j=0}^i a_1(j)\left (\sum_{k=0}^{n-i-1}
  \underline{\mathfrak{S}}^{(n-2)}_{k,i-1} x_n^{2k}\right )+\sum_{j=0}^i a_1(j)\left (\sum_{k=0}^{n-i-1}
  \underline{\mathfrak{S}}^{(n-2)}_{k,i-1} x_n^{2k}\right )+\\x_1\left
  (\sum_{j=1}^{2i-1}a_2(j) x_2\cdots x_{n-1} x_n^k +\sum_{j=0}^i
  a_1(j)\sum_{k=0}^{n-i}\underline{\mathfrak{S}}^{(n-2)}_{k,i-2}x_n^{2k}\right)
\end{multline*}
then $G^{(n)}_{i+1}=f_{i+1}\left ( \sum_{k=0}^i a_1(k)
  x^{2k}_n\right)-G^{(n)}_i x_n$.
\end{lem}

\begin{proof}  Let $h_1:=S(f_{i+1},G^{(n)}_i)$.  Recall that the
  degree of $f_{i+1}$ is
  $(\underbrace{2,\dots,2}_{n-i},\underbrace{1,\dots,1}_{i})$ and that
  the leading coefficient of $f_{i+1}$ is 1.  So, $h_1=x_n^{2i-1}a_1(i)
  f_{i+1}-G^{(n)}_i$.  Hence, $deg(h_1)<
  (\underbrace{2,\dots,2}_{n-i},\underbrace{1,\dots,1}_{i-1},2i)$ and
  $deg(h_1)\geq
  (\underbrace{2,\dots,2}_{n-i},\underbrace{1,\dots,1}_{i-1},2i-2)$
  since a term with this exponent appears in $G^{(n)}_i$ but not in $x_n^{2i-1}a_1(i)
  f_{i+1}$ (if it did, then there would have to be a term in
  $f_{i+1}$ with a multiplicand of $x_n^{-1}$).  Note that
  $deg(h_1)\neq
  (\underbrace{2,\dots,2}_{n-i},\underbrace{1,\dots,1}_{i-1},2i-1)$
  because if it did, there would be a summand of $f_{i+1}$ with
  degree
  $(\underbrace{2,\dots,2}_{n-i+1},\epsilon_1,\dots,\epsilon_{i-1})>deg(f_{i+1})$.

For $k>1$ assume that 
\[
h_{k-1}=f_{i+1}(x_n^{2i-1}a_1(i)+\cdots+x_n^{2i-2k+3}a_1(i-k+2))-G^{(n)}_i.
\]
Now, $h_k=S(f_{i+1},h_{k-1})$.  Notice that
\[
(\underbrace{2,\dots,2}_{n-i},\underbrace{1,\dots,1}_{i-1},2i-2k)\leq
deg(h_k)<(\underbrace{2,\dots,2}_{n-i},\underbrace{1,\dots,1}_{i-1},2i-2k+2).
\]
We deduce the strict inequality from the definition of the
$S$-polynomial and the other inequality from the fact that a term of
degree $(\underbrace{2,\dots,2}_{n-i},\underbrace{1,\dots,1}_{i-1},2i-2k)$ appears in $G^{(n)}_i$.  Moreover, we cannot have a term of
degree
$(\underbrace{2,\dots,2}_{n-i},\underbrace{1,\dots,1}_{i-1},2i-2k+1)$
for that would imply there was a term in $f_{i+1}$ of degree
$(\underbrace{2,\dots,2}_{n-i+1},\epsilon_1,\dots,\epsilon_{i-1})>deg(f_{i+1})$.
Thus $deg(h_k)=(\underbrace{2,\dots,2}_{n-i},\underbrace{1,\dots,1}_{i-1},2i-2k)$.

Notice when $i=k$, $h_i$ has degree
$(\underbrace{2,\dots,2}_{n-i},\underbrace{1,\dots,1}_{i-1},0)$ and
$h_i=\sum_{k=1}^i x_n^{2k-1}a_1(k) f_{i+1}-G^{(n)}_i$.  Then
\begin{align*}
G^{(n)}_{i+1}=h_{i+1}&=a_1(0) f_{i+1}+x_nh_i\\
&=a_1(0)f_{i+1}+x_n(\sum_{k=1}^i x_n^{2k-1}a_1(k) f_{i+1}-G^{(n)}_i)\\
&=f_{i+1}\left ( \sum_{k=0}^i a_1(k)x_n^{2k}\right)-G^{(n)}_ix_n.
\end{align*}\end{proof}

\begin{lem}\label{lem3}  Let $G^{(n)}_i$ be as in the statement of
  \thmref{algorithm}.  Then
\begin{multline*}
G^{(n)}_i=x_1^2\sum_{j=0}^i a_1(j)\left (\sum_{k=0}^{n-i-1}
  \underline{\mathfrak{S}}^{(n-2)}_{k,i-1} x_{n}^{2k}\right )+\sum_{j=0}^i a_1(j)\left (\sum_{k=0}^{n-i-1}
  \underline{\mathfrak{S}}^{(n-2)}_{k,i-1} x_{n}^{2k}\right )\\+x_1\left
  (\sum_{j=1}^{2i-1}a_2(j) x_2\cdots x_{n-1} x_n^k +\sum_{j=0}^i
  a_1(j)\sum_{k=0}^{n-i}\underline{\mathfrak{S}}^{(n-2)}_{k,i-2}x_{n}^{2k}\right)
\end{multline*}
\noindent where the sequences $a_1(k)$ and $a_2(k)$ are palindromic.
\end{lem}
\begin{proof}
We prove this lemma by induction on $i$ and for $n > 3$ (the examples
in \secref{alg-ex} show the lemma for $n=2$ and $n=3$).  We assume
that $G^{(n)}_i$ is in the form prescribed by the lemma and then use
\lemref{lem2} to show that $G^{(n)}_{i+1}$ is of the same form.

In \lemref{lem1} we broke up $f_{i+1}$ according to the power of $x_n$
in the terms, and we rewrite it as in equation \ref{f-i-4}
\begin{multline*}\label{f-i-5}
f_{i+1}:=(x_1^2+1)\left( (x_{n}^2+1)\sum_{k=0}^{n-i-2}
  \underline{\mathfrak{S}}^{(n-2)}_{k,i}+x_n\sum_{k=0}^{n-i-1}\underline{\mathfrak{S}}^{(n-2)}_{k,i-1}\right)+\\
x_n\left((x_{n}^2+1)\sum_{k=0}^{n-i-1}\underline{\mathfrak{S}}^{(n-2)}_{k,i-1}+x_{n}(\sum_{k=0}^{n-i}\underline{\mathfrak{S}}^{(n-2)}_{k,i-2}-c_{i+1}x_2\cdots
  x_{n-1})\right).
\end{multline*}

By \lemref{lem2} we know
\begin{equation}\label{G}
G^{(n)}_{i+1}=f_{i+1}\left ( \sum_{k=0}^i a_1(k)
  x_n^{2k}\right)-G^{(n)}_i x_n.
\end{equation}
We look at the $x_1^2$ term in $G^{(n)}_{i+1}$:
\begin{align*}
x_1^2&\left ( \sum_{k=0}^i a_1(k) x_n^{2k}\right)\left ( (x_{n}^2+1)\sum_{k=0}^{n-i-2}
  \underline{\mathfrak{S}}^{(n-2)}_{k,i}+x_n\sum_{k=0}^{n-i-1}\underline{\mathfrak{S}}^{(n-2)}_{k,i-1}\right
  )-\\
&x_n x_1^2\left((\sum_{k=0}^i a_1(k)
  x_n^{2k})(\sum_{k=0}^{n-i-1}\underline{\mathfrak{S}}^{(n-2)}_{k,i-1})\right)\\
&=x_1^2\left((\sum_{k=0}^i a_1(k)x_n^{2k+2})(\sum_{k=0}^{n-i-2}
  \underline{\mathfrak{S}}^{(n-2)}_{k,i})+(\sum_{k=0}^i a_1(k)x_n^{2k})(\sum_{k=0}^{n-i-2}
  \underline{\mathfrak{S}}^{(n-2)}_{k,i})\right)\\
&:=x_1^2\left (\sum_{k=0}^{i+1} b_1(k)x_n^{2k} (\sum_{k=0}^{n-i-2}
  \underline{\mathfrak{S}}^{(n-2)}_{k,i})\right)
\end{align*}
\noindent where $b_1(0)=a_1(0)$, $b_1(i+1)=a_1(i)$ and for $ 0<\ell<i+1$,
$b_1(\ell)=a_1(\ell)+a_1(\ell-1)$.  We need to show that the sequence
$b_1(k)$ is palindromic; i.e., we want to show $b_1(k)=b_1(i+1-k)$.
We check the endpoints of the sequence first:  since $a_1(k)$ is palindromic, we know that
$b_1(0)=a_1(0)=a_1(i)=b_1(i+1)$.  Now let $ 0<\ell<i+1$.  Then
\begin{align*}
b_1(\ell)&=a_1(\ell)+a_1(\ell-1)\\
&=a_1(i-\ell)+a_1(i-\ell+1)\\
&=b_1(i-\ell+1).
\end{align*}
\noindent  Thus the $x_1^2$ term of equation \ref{G} is of the form
demanded by the lemma.

Since the $x_1^2$ terms in both $x_n G^{(n)}_i$ and
$f_{i+1}\sum_{k=0}^i a_1(k)x_n^{2k}$ are identical to the $x_1^0$
terms in  $x_n G^{(n)}_i$ and
$f_{i+1}\sum_{k=0}^i a_1(k)x_n^{2k}$, the same argument as above shows
that the $x_1^0$ term of equation \ref{G} is of the form demanded by
the lemma.

Now we turn our attention to the $x_1$ term in equation \ref{G}.  We
have
\begin{align*}
x_1&\left ( \sum_{k=0}^i a_1(k) x_n^{2k}\right)\left ( (x_n^2+1)(\sum_{k=0}^{n-i-1}\underline{\mathfrak{S}}^{(n-2)}_{k,i-1})+x_{n}(\sum_{k=0}^{n-i}\underline{\mathfrak{S}}^{(n-2)}_{k,i-2}-c_{i+1}x_2\cdots
  x_{n-1})\right)-\\
&x_1 x_n \left( \sum_{k=1}^{2i-1} a_2(k) x_2\cdots x_{n-1}
  x_n^{k}+(\sum_{k=0}^i
  a_1(k)x_n^{2k})(\sum_{k=0}^{n-i}\underline{\mathfrak{S}}^{(n-2)}_{k,i-2})\right)\\
&=x_1( (\sum_{k=0}^i a_1(k)
  x_n^{2k+2})(\sum_{k=0}^{n-i-1}\underline{\mathfrak{S}}^{(n-2)}_{k,i-1})+
  (\sum_{k=0}^i a_1(k)
  x_n^{2k})(\sum_{k=0}^{n-i-1}\underline{\mathfrak{S}}^{(n-2)}_{k,i-1})+\\&\sum_{k=1}^{2i-1}
  a_2(k)x_2\cdots x_{n-1}x_n^{k+1}-c_{i+1}\sum_{k=0}^i a_1(k)
  x_2\cdots x_{n-1} x_n^{2k+1})\\
&:=x_1\left((\sum_{k=0}^{i+1} b_1(k) 
  x_n^{2k})(\sum_{k=0}^{n-i-1}\underline{\mathfrak{S}}^{(n-2)}_{k,i-1})+(\sum_{k=1}^{2i+1}
  b_2(k)x_2\cdots x_{n-1}x_n^{k})\right)
\end{align*}
\noindent where $b_1(k)$ was defined above and $b_2(\ell)$ is $a_2(\ell-1)$
  if $\ell$ is even and $a_2(\ell-1)-c_{i+1}
  a_1\left(\frac{\ell-1}{2}\right)$ if $\ell$ is odd.  We have already
  shown that $b_1(k)$ is palindromic so all that remains to be shown
  is that $b_2(\ell)$ is palindromic. Since $a_2(k)$ and $a_1(k)$ are
  palindromic (i.e., $a_2(k)=a_2(2i-k)$ and $a_1(k)=a_1(i-k)$ we get,
  for odd $\ell$, 
\begin{align*}
b_2(\ell)&=a_2(\ell-1)-c_{i+1} a_1\left(\frac{\ell-1}{2}\right)\\
&=a_2(2i-\ell+1)-c_{i+1}a_1\left(\frac{2i-k+1}{2}\right)\\
&=a_2(2i-\ell+2-1)-c_{i+1}a_1\left(\frac{2i-\ell+2-1}{2}\right)\\
&=b_2(2i-\ell+2);
\end{align*}
and for even $\ell$ we get
\begin{align*}
b_2(\ell)&=a_2(\ell-1)\\
&=a_2(2i-\ell-k+2-1)\\
&=b_2(2i-\ell+2).
\end{align*}
\noindent Thus, the $x_1$ term of equation \ref{G} is as described in
the lemma and the lemma is proved.  
\end{proof}

\begin{proof}[Proof~of~\thmref{algorithm}]
By \lemref{lem3} 
\begin{align*}
G^{(n)}_n=x_1^2&\sum_{j=0}^n a_1(j)\left (\sum_{k=0}^{-1}
  \underline{\mathfrak{S}}^{(n-2)}_{k,n-1} x_{n}^{2k}\right )+\sum_{j=0}^n a_1(j)\left (\sum_{k=0}^{-1}
  \underline{\mathfrak{S}}^{(n-2)}_{k,n-1} x_{n}^{2k}\right )+\\
&x_1\left
  (\sum_{j=1}^{2n-1}a_2(j) x_2\cdots x_{n-1} x_n^k +\sum_{j=0}^i
  a_1(j)\sum_{k=0}^{n-i}\underline{\mathfrak{S}}^{(n-2)}_{k,n-2}x_{n}^{2k}\right)\\
&=\sum_{k=1}^{2n-1} a_2(k) x_1\cdots x_{n-1}x_n^k+\sum_{k=0}^n a_1(k)
  x_1\cdots x_{n-1}x_n^{2k}\\
&:=\sum_{k=0}^{2n} a(k) x_1\cdots x_{n-1}x_n^k
\end{align*}
\noindent where $a_1(k)$ and $a_2(k)$ are palindromic and
  $a(0)=a_1(0)$, $a(2n)=a_n(0)=a_1(0)$, for even
  $0<\ell<2n$, $a(\ell)=a_1\left(\frac{\ell}{2}\right)+a_2(\ell)$ and
  for odd $0<\ell<2n$, $a(\ell)=a_2(\ell)$.  For even $0<\ell <2n$ we
  have
\begin{align*}
a(\ell)&=a_1\left(\frac{\ell}{2}\right)+a_2(\ell)\\
&=a_1\left(n-\frac{\ell}{2}\right)+a_2(2n-\ell)\\
&=a_2(2n-\ell)+a_1\left(\frac{2n-\ell}{2}\right)\\
&=a(2n-\ell)
\end{align*}
and for odd $0<\ell<2n$ we have
\begin{align*}
a(\ell)&=a_2(\ell)\\
&=a_2(2n-\ell)\\
&=a(2n-\ell).
\end{align*}
\noindent  This proves the theorem.
\end{proof}

Immediately, we deduce the following corollary:
\begin{cor}\label{Pn-palindromic}
Let the notation be as in \thmref{algorithm}.  The polynomial $P_n:=G^{(n)}_n/x_1\cdots x_{n-1}$ is palindromic.
\end{cor}

By the following proposition, being able to find $P_n$ allows us to
find all the solutions to system \ref{sys4}:
\begin{prop}\label{num-solns}  There are $2^n n!$ solutions
  $(\alpha_1,\dots,\alpha_n)\in(\bbC^\times)^n$ to system \ref{sys4} produced
  by algorithm $\mathcal{A}$.
\end{prop}
\begin{proof}  The proof is by induction.  When $n=1$ we see that
  $f_1=G_1^{(1)}$ is a quadratic polynomial in $x_1$ and thus there are 2 solutions. 

Let $n>1$ and assume the proposition holds for $m<n$.  Apply algorithm
$\mathcal{A}$ to system \ref{sys4} to derive the palindromic
polynomial $P_n$ of degree $2n$ with roots
$z_1^{(n)},\dots,z_{2n}^{(2n)}$.  Without loss of generality, we substitute $z_1^{(n)}$ into system \ref{sys4} and, rewriting each $f_i$ according
to equation \ref{f-i-1}, we derive the following associated system:
\begin{align}\label{hat-sys0}
&\hat{f}_1:=((z_1^{(n)})^2+1)\sum_{k=0}^{n-1} \mathfrak{S}^{(n-1)}_{k,0}-z_1^{(n)}c_1x_1\dots x_{n-1} =0\nonumber\\
&\hat{f}_2:=((z_1^{(n)})^2+1)\sum_{k=0}^{n-2}\mathfrak{S}^{(n-1)}_{k,1}+z_1^{(n)}(\sum_{k=0}^{n-1}\mathfrak{S}^{(n-1)}_{1,0}- c_2x_1\dots x_{n-1})=0\nonumber\\
&\vdots\nonumber\\
&\hat{f}_i:=((z_1^{(n)})^2+1)\sum_{k=0}^{n-i}\mathfrak{S}^{(n-1)}_{k,i-1}+z_1^{(n)}(\sum_{k=0}^{n-i+1}\mathfrak{S}^{(n-1)}_{k,i-2}- c_ix_1\dots x_{n-1})=0\\
&\vdots\nonumber\\
&\hat{f}_n:=((z_1^{(n)})^2+1)\mathfrak{S}^{(n-1)}_{k,n-1}+z_1^{(n)}(\sum_{k=0}^{1}\mathfrak{S}^{(n-1)}_{k,n-2}- c_nx_1\dots x_{n-1})=0\nonumber .
\end{align}

Define $\hat{f}_{n+1}:=z_1^{(n)}\hat{f}_1-((z_1^{(n)})^2+1)\hat{f}_2$
and for $i\geq 2$,
\[
\hat{f}_{n+i}=\hat{f}_{i+1}((z_1^{(n)})^2+1)^i-\hat{f}_{n+i-1}z_1^{(n)}.
\]
We claim
$\hat{f}_{n+i}=((z_1^{(n)})^2+1)^{i+1}\sum_{k=0}^{n-i-1}\mathfrak{S}^{(n-1)}_{k-1,i}-\hat{c}_i
x_1\cdots
x_{n-1}$.
\begin{proof}[Proof of Claim]  The claim is proved by induction on
  $i$.  By definition and system \ref{hat-sys0} we have 
\begin{align*}
\hat{f}_{n+1}&= z_1^{(n)}\hat{f}_1-((z_1^{(n)})^2+1)\hat{f}_2\\
&=z_1^{(n)}\left(((z_1^{(n)})^2+1)\sum_{k=0}^{n-1}
\mathfrak{S}^{(n-1)}_{k,0}-z_1^{(n)}x_1\dots
x_{n-1}\right)-\\
&((z_1^{(n)})^2+1)\left(((z_1^{(n)})^2+1)\sum_{k=0}^{n-2}\mathfrak{S}^{(n-1)}_{k,1}+z_1^{(n)}(\sum_{k=0}^{n-1}\mathfrak{S}^{(n-1)}_{k,0}-
  c_2x_1\dots x_{n-1}) \right)\\
&=((z_1^{(n)})^2+1)^2\sum_{k=0}^{n-2}\mathfrak{S}^{(n-1)}_{k,1}-((z_1^{(n)})^2+1)(z_1^{(n)})c_2x_1\cdots
x_{n-1}-\\&(z_1^{(n)})^2 c_1\cdots x_1\cdots x_{n-1}
\end{align*}
Setting $\hat{c}_1:=(z_1^{(n)})^2 c_1-((z_1^{(n)})^2+1)(z_1^{(n)})c_2$
we can write 
\[
\hat{f}_{n+1}=((z_1^{(n)})^2+1)^2\sum_{k=0}^{n-2}\mathfrak{S}^{(n-1)}_{k,1}-\hat{c}_1
x_1\cdots x_{n-1}
\]
as demanded by the claim.  
\end{proof}
Now assume that 
\[
\hat{f}_{n+i-1}=((z_1^{(n)})^2+1)^i\sum_{k=0}^{n-i}\mathfrak{S}^{(n-1)}_{k,i-1}-\hat{c}_{i-1}
x_1\cdots x_{n-1}.
\]
By definition and an argument similar to the case when $i=1$ we get
\begin{align*}
\hat{f}_{n+i}&=((z_1^{(n)})^2+1)^i\hat{f}_{i+1}-(z_1^{(n)})\hat{f}_{n+i-1}\\
&=((z_1^{(n)})^2+1)^i\left(((z_1^{(n)})^2+1)\sum_{k=0}^{n-i-1}\mathfrak{S}^{(n-1)}_{k,i}+z_1^{(n)}(\sum_{k=0}^{n-i}\mathfrak{S}^{(n-1)}_{k,i-2}-
  c_{i+1}x_1\dots x_{n-1})\right)-\\
&(z_1^{(n)})\hat{f}_{n+i-1}\\
&=((z_1^{(n)})^2+1)^i\left(((z_1^{(n)})^2+1)\sum_{k=0}^{n-i-1}\mathfrak{S}^{(n-1)}_{1,i}+z_1^{(n)}(\sum_{k=0}^{n-i}\mathfrak{S}^{(n-1)}_{k,i-2}-
  c_{i+1}x_1\dots x_{n-1})\right)-\\
&(z_1^{(n)}\left(((z_1^{(n)})^2+1)^i\sum_{k=0}^{n-i}\mathfrak{S}^{(n-1)}_{k,i-1}-\hat{c}_{i-1}
x_1\cdots x_{n-1}\right)\\
&=((z_1^{(n)})^2+1)^{i+1}\left(\sum_{k=0}^{n-i-1}\mathfrak{S}^{(n-1)}_{k,i}\right)-((z_1^{(n)})^2+1)^i(z_1^{(n)})c_{i+1}x_1\dots
x_{n-1}-\\
&(z_1^{(n)})c_{i-1}x_1\dots x_{n-1}
\end{align*}
\noindent which is of the desired form.
Note that the system
\begin{align}\label{hat-sys1}
&\frac{\hat{f}_{n+1}}{((z_1^{(n)})^2+1)^2}=\sum_{k=0}^{n-2}
\mathfrak{S}^{(n-1)}_{k-1,1}-\frac{\hat{c}_1}{((z_1^{(n)})^2+1)^2}x_1\cdots
x_{n-1}\nonumber\\
&\vdots\nonumber\\
&\frac{\hat{f}_{n+i}}{((z_1^{(n)})^2+1)^{i+1}}=\sum_{k=0}^{n-i-1}
\mathfrak{S}^{(n-1)}_{k-1,i}-\frac{\hat{c}_i}{((z_1^{(n)})^2+1)^{i+1}}x_1\cdots
x_{n-1}\\
&\vdots\nonumber\\
&\frac{\hat{f}_{2n-1}}{((z_1^{(n)})^2+1)^{n}}=
\mathfrak{S}^{(n-1)}_{0,n-1}-\frac{\hat{c}_{n-1}}{((z_1^{(n)})^2+1)^{n}}x_1\cdots
x_{n-1}\nonumber
\end{align}
\noindent is of the same form has the input for algorithm
$\mathcal{A}$.  Thus, by induction, algorithm $\mathcal{A}$ produces
$2^{n-1}(n-1)!$ $(n-1)$-tuples that solve system \ref{hat-sys1}.  Since
there are $2n$ possible choices for $z_i^{(n)}$ we get a total of $2^n
n!$ solutions to system \ref{sys4}.
\end{proof}
\begin{cor}\label{Sat-comp-cor}  If system \ref{sys4} arises from a simultaneous eigen
  cusp form $F$, algorithm $\mathcal{A}$ produces the Satake
  $p$-parameters of $F$.
\end{cor}
\begin{proof}
According to \propref{num-solns}, algorithm $\mathcal{A}$ generates
$2^n n!$ $n$-tuples in $(\bbC^\times)^n$ that solve the equations in system
\ref{sys4}.  Pick one such $n$-tuple, call it
$(\alpha_1,\dots,\alpha_n)$.  Then, as in system \ref{sys1},
we get
\[
\frac{\lambda_{F,0}(p)}{(1+\alpha_1)\cdots (1+\alpha_n)}=\alpha_0.
\]
Next, since the second equation in \ref{sys1},
$\lambda_n=c_{11}\alpha_0^2\alpha_1\cdots\alpha_n$, is
nonzero on the left hand side, each $\alpha_i$ is nonzero.  Finally,
since the equations in system \ref{sys1} are $W_n$ invariant, we
conclude that $(\alpha_0,\dots, \alpha_n)\in(\bbC^\times)^{n+1}/W_n$.
Moreover, there is only one such $(n+1)$-tuple since we found $2^n n!$
solutions and the cardinality of $W_n$ is $2^n n!$.
\end{proof}

\section{Example of Implementing the Algorithm}\label{alg-ex}

We conclude by providing some examples.  The first example
illustrates the implementation of the algorithm for $n=3$ and the
second applies the implementation of the algorithm for $n=2$ and
applies the results of the algorithm to the verification of the
Generalized Ramanujan Petersson Conjecture which states the Siegel
modular forms not in the Maass Space have unimodular Satake
$p$-parameters.  We also provide a table of Satake $p$-parameters for
the Schottky form $J\in\cS_8(\Gamma_4)$ \cite{Breul2} for small $p$.
In presenting these examples we also prove the base case for the proof
of \lemref{lem3}.

\subsection{Sample Implementation}
Suppose $n=3$ and recall that the matrix
$[\Omega_2]=(c_{ij})$ is upper triangular.  By \propref{Aprop} we know the eigenvalues of
$F\in\cS_k(\Gamma_3)$ with respect to the operators $T_0(p)$,
$T_0(p^2)$, $T_1(p^2)$, $T_2(p^2)$ and $T_3(p^2)$ are polynomials in the Satake
parameters:
\begin{align}\label{sys1-n=3}
\lambda_{F,0}(p)&=x_0(1+x_1)(1+x_2)(1+x_3)\nonumber\\
\lambda_3&=c_{11}[x_0^2x^{(1,1,1)}]\nonumber\\
\lambda_2&=c_{12}[x_0^2x^{(1,1,1)}]+c_{22}[x_0^2x^{(2,1,1)}]\\
\lambda_1&=c_{13}[x_0^2x^{(1,1,1)}]+c_{23}[x_0^2x^{(2,1,1)}]+c_{33}[x_0^2x^{(2,2,1)}]\nonumber\\
\lambda_0&=c_{14}[x_0^2x^{(1,1,1)}]+c_{24}[x_0^2x^{(2,1,1)}]+c_{34}[x_0^2x^{(2,2,1)}]+c_{44}[x_0^2x^{(2,2,2)}]\nonumber
\end{align}
The first equation in system \ref{sys1-n=3} will be used to compute
the value of $x_0$ and we ignore it for the time being.

We define constants $k_1$, $k_2$, $k_3$ and $k_4$ to be the following:
\begin{align}\label{sys2-n=3}
\lambda_3&=c_{11}[x_0^2x^{(1,1,1)}]\nonumber\\&:=k_1c_{11}\nonumber\\
\lambda_2&=c_{12}[(x_0^2x^{(1,1,1)}]+c_{22}[x_0^2x^{(2,1,1)}]\nonumber\\&:=c_{12}k_1+c_{22}k_2\nonumber\\
\lambda_1&=c_{13}[x_0^2x^{(1,1,1)}]+c_{23}[x_0^2x^{(2,1,1)}]+c_{33}[x_0^2x^{(2,2,1)}]\\&:=k_1c_{13}+k_2c_{23}+k_3c_{33}\nonumber\\
\lambda_0&=c_{14}[x_0^2x^{(1,1,1)}]+c_{24}[x_0^2x^{(2,1,1)}]+c_{34}[x_0^2x^{(2,2,1)}]+c_{44}[x_0^2x^{(2,2,2)}]
\nonumber\\&:=k_1c_{14}+k_2c_{24}+k_3c_{34}+k_4c_{44}\nonumber
\end{align}

Now we rewrite the equations in system \ref{sys2-n=3} as
\begin{align}\label{sys3-n=3}
&k_1=\frac{\lambda_{3}}{c_{11}}=[x_0^2x^{(1,1,1)}]=x_0^2x_1x_2x_3\nonumber\\
&k_2=\frac{\lambda_{2}-c_{12}k_1}{c_{22}}=x_0^2\sum_{k=0}^1\mathfrak{S}^{(3)}_{k,2}\\
&k_3=\frac{\lambda_{1}-c_{13}k_1-c_{23}k_2}{c_{33}}=x_0^2\sum_{k=0}^2\mathfrak{S}^{(3)}_{k,1}\nonumber\\
&k_4=\frac{\lambda_{0}-c_{14}k_1-c_{24}k_2-c_{34}k_3}{c_{44}}=x_0^2\sum_{k=0}^3\mathfrak{S}^{(3)}_{k,0}\nonumber
\end{align}

Finally to eliminate the variable $x_0$ we rewrite the equations in
\ref{sys3-n=3} and drop the first equation:
\begin{align}\label{sys4-n=3}
f_1:=\sum_{k=0}^3\mathfrak{S}^{(3)}_{k,0}-c_1x_1x_2x_3=0\nonumber\\
f_2:=\sum_{k=0}^3\mathfrak{S}^{(3)}_{k,1}-c_2x_1x_2x_3=0\\
f_3:=\sum_{k=0}^3\mathfrak{S}^{(3)}_{k,2}-c_1x_1x_2x_3=0\nonumber.
\end{align}

As described in \thmref{algorithm} Algorithm $\mathcal{A}$ produces
polynomials $f_4$, $f_5$, $f_6$, $f_7$ and $f_8$, and a final
polynomial that is of degree $(1,1,6)$, palindromic, and divisible by
$x_1x_2$.  These calculations were done with Maple 9 and checked by hand:
\begin{align*}
f_4:=&S(f_1,f_2)\\
=&x^{(2,2,0)}-x^{(2,1,3)}-x^{(2,1,1)}+x^{(2,0,0)}-x^{(1,2,3)}-x^{(1,2,1)}-c_1x^{(1,1,1)}+c_2x^{(1,1,2)}\\&-x^{(1,0,3)}-x^{(1,0,1)}+x^{(0,0,0)}+x^{(0,2,0)} -x^{(0,1,1)}-x^{(0,1,3)}\\
f_5:=&S(f_2,f_4)\\
=&x^{(2,1,4)}+2x^{(2,1,2)}+x^{2,1,0}+x^{(1,2,4)}+2x^{(1,2,2)}+x^{(1,2,0)}-c_2x^{(1,1,3)}+c_1x^{(1,1,2)}\\&-c_2x^{(1,1,1)}+x^{(1,0,4)}+2x^{(1,0,2)}+x^{(1,0,0)}+x^{(0,1,4)}+2x^{(0,1,2)}+x^{(0,1,0)}\\
f_6:=&S(f_3,f_5)\\
=&2x^{(2,1,2)}+x^{2,1,0)}+x^{(1,2,0)}+2x^{(1,2,2)}-x^{(1,1,5)}+c_3x^{(1,1,4)}+(-c_2-1)x^{(1,1,3)}\\&+c_1x^{(1,1,2)}-c_2x^{(1,1,1)}+x^{(1,0,0)}+2x^{(1,0,2)}\\&+2x^{(0,1,2)}+x^{(0,1,0)}\\
f_7:=&S(f_3,f_6)\\
=&x^{(2,1,0)}+x^{(1,2,0)}-x^{(1,1,5)}+c_3x^{(1,1,4)}+(-c_2-3)x^{(1,1,3)}+(2c_3+c_1)x^{(1,1,2)}\\&+(-c_2-2)x^{(1,1,1)}+x^{(1,1,0)}+y^{(0,1,0)}\\
f_8:=&S(f_3,f_7)\\
=&-x^{(1,1,6)}+c_3x^{(1,1,5)}+(-c_2-3)x^{(1,1,4)}+(2c_3+c_1)x^{(1,1,3)}\\&+(-c_2-3)x^{(1,1,2)}+c_3x^{(1,1,1)}-x^{(1,1,0)}.
\end{align*} 
Thus 
\begin{align*}
P_3&:=G^{(3)}_3/x_1x_2=f_8/x_1x_2\\
&=-x_3^6+c_3x_3^5+(-c_2-3)x_3^4+(2c_3+c_1)x_3^3+(-c_2-3)x_3^2+c_3x_3-1
\end{align*}
\noindent is palindromic as demanded by \corref{Pn-palindromic}.

Suppose the roots of $P_3\in\bbC[x_3]$ are $\gamma_i$ for $1\leq i \leq 6$.  Then as in
the proof of \propref{num-solns} we turn system \ref{sys4-n=3} into
the following system of equations:
\begin{align}\label{sys5-n=3}
\hat{f}_1&:=(\gamma_1^2+1)\sum_{k=0}^2\mathfrak{S}^{(2)}_{k-1,0}-c_1\gamma_1x_2=0\nonumber\\
\hat{f}_2&:=(\gamma_1^2+1)\sum_{k=0}^1\mathfrak{S}^{(2)}_{k-1,1}+\gamma_1(\sum_{k=0}^2\mathfrak{S}^{(2)}_{k-1,0}-c_2x_1x_2)=0\nonumber\\
\hat{f}_3&:=(\gamma_1^2+1)x_1x_2+\gamma_1(\sum_{k=0}^1\mathfrak{S}^{(2)}_{k-1,1}-c_3x_1x_2)=0\nonumber\\
\hat{f}_4&:=\gamma_1\hat{f}_1-(\gamma_1^2+1)\hat{f}_2\\
&= (\gamma_1^2+1)^{2}\sum_{k=0}^{1}\mathfrak{S}^{(2)}_{k-1,1}-\hat{c}_1
x_1\cdots x_{2}=0\nonumber\\
\hat{f}_5&:=\hat{f}_{2}(\gamma_1^2+1)^2-\hat{f}_{4}\gamma_1\nonumber\\
&=(\gamma_1^2+1)^{3}\mathfrak{S}^{(2)}_{0,2}-\hat{c}_2
x_1\cdots x_2=0\nonumber
\end{align}
\noindent for complex constants $\hat{c}_1$ and $\hat{c}_2$.  

Now we can apply algorithm $\mathcal{A}$ to the last two equations in
system \ref{sys5-n=3} to get, as in \corref{Pn-palindromic}, a palindromic polynomial
$P_2\in\bbC[x_2]$ of degree 4.  Suppose $\beta_i$ for $1\leq i\leq 4$ are
the roots of $P_2$.  Next, we choose, without loss of generality, the
root $\beta_1$ and substitute it into the first two equations, and produce
\[
\hat{\hat{f}}_1:=(\gamma_1^2+1)(\beta_1^2+1)(x_1^2+1)-c_1\beta_1\gamma_1x_1=0,
\]
a palindromic quadratic polynomial in $\bbC[x_1]$.  Finally, suppose
that $\alpha_1$ and $\alpha_2$ are the roots of $\hat{\hat{f}}_1$.  

Recall that in this case, the Satake $p$-parameters are an element of
$(\bbC^\times)^4/W_3$, so we have yet to find the
$p$-parameter $\alpha_{0,p}$.  From the first equation of \ref{sys1} we get
\[
\lambda_{F,0}(p)=x_0(1+x_1)(1+x_2)(1+x_3).
\]
If we make the substitutions $x_1=\alpha_1$, $x_2=\beta_1$ and
$x_3=\gamma_1$ we can solve for the first $p$-parameter by
\[
\delta=\frac{\lambda_{F,0}(p)}{(1+\alpha_1)(1+\beta_1)(1+\gamma_1)}
\]
and conclude from \corref{Sat-comp-cor} that
$(\delta,\alpha_1,\beta_1,\gamma_1)$ are the Satake $p$-parameters.

Since we have not included a computation of the matrix $[\Omega_2]$
for $n=3$ we will not compute the Satake $p$-parameters of Hecke eigenforms $\cS_k(\Gamma_3)$.  In the following section we
compute the Satake $p$-parameters for particular Hecke eigenforms of
genus 2 and 4.

\subsection{Satake $p$-parameters for genus $2$}

In this section we compute, for small primes $p$, the Satake
$p$-parameters of particular cusp forms not in the Maass space.  We
define the so-called Maass space:
\begin{definition}
 Let $k\in\bbZ$ be even.  The \textit{Maass space}
  $S_k^*(\Gamma_2)$ is
\[
\left\{ F\in S_k(\Gamma_2):F=\iota_k(f) \text{ for some $f\in
    S_k(\Gamma_1)$}\right\}
\]
where $\iota_k$ is the Saito-Kurokawa lift.
\end{definition}

Skoruppa (\cite{Skoruppa}) was able to compute the Hecke eigenvalues
for the first 6 cusp forms of even weight which do not belong the
Maass space.  The first of the forms occurs at weight 20 and he calls
it $\Upsilon 20$.  Similarly, there is one such cusp form of weight 22,
$\Upsilon 22$; two of weight 24, $\Upsilon 24a$ and $\Upsilon 24b$;
and two of weight 26, $\Upsilon 26a$ and $\Upsilon 26b$.  By
computing the Fourier coefficients of these Hecke eigenforms, he
was able to generate the Table \ref{evalues} of eigenvalues.  

From applying algorithm $\mathcal{A}$, we were able to generate Table \ref{tab-Sat-par} of Satake $p$-parameters.  Note, in particular, that
all the Satake $p$-parameters have modulus 1.  A piece of notation
needed to interpret the table is the following Hecke operator: 
\[
T(p^2)=\sum_{i=0}^n T_i(p^2).
\]

\begin{table}
\begin{tabular}{|l|r|r|r|}
\hline
$\Upsilon *$ & $p$ & $\lambda_{\Upsilon *}(T(p))$ &
$\lambda_{\Upsilon *}(T(p^2))$\\
\hline
$\Upsilon 20$ & $2$ & $-2^8\cdot 3^2\cdot 5\cdot 73$ & $2^{16} \cdot
523\cdot 7243$\\
& $3$ & $2^3\cdot 3^5\cdot 5\cdot 7\cdot 5099$ & $-3^{10} \cdot
2658457\cdot 2879687$\\
& $5$ & $-2^2\cdot 3^2\cdot 5^3\cdot 7\cdot 166103087$& $-5^6\cdot 9973\cdot
1165906151989603 $\\
& $7$ & $2^4\cdot 5^2\cdot 7^3\cdot 673\cdot 28346749 $ & $-3\cdot 7^6
\cdot 23\cdot 43\cdot 71\cdot 1327$\\& & &$\cdot 1844737
\cdot 90682160593$\\
\hline
$\Upsilon 22$ & $2$ & $-2^8\cdot 3\cdot 5\cdot 577$& $2^{16}\cdot
18869089$\\
& $3$ & $-2^3\cdot 3^5\cdot 5\cdot 19 \cdot 97\cdot 167$ &
$3^{10}\cdot 797\cdot 1049\cdot 2707 \cdot 48271 $\\
& $5$ & $2^2\cdot 3 \cdot 5^3\cdot 60700091989$& $-5^6\cdot 105263022561216721922951$\\
\hline
$\Upsilon 24a$&$2$&$-2^{11}\cdot 3\cdot 5\cdot 181$& $2^{22}\cdot
1039321$\\
&$3$&$-2^3\cdot 3^6\cdot 5\cdot 7\cdot 23^2\cdot 491$&$3^{12}\cdot
7\cdot 9027177753487$\\
& $5$ &$-2^2 \cdot 3\cdot 5^3\cdot 7\cdot 29\cdot
109438961$&$-2^2\cdot 3\cdot 5^3\cdot 7\cdot 29 \cdot 109438961$\\
\hline
$\Upsilon 24b$&$2$&$-2^9 \cdot 3^2 \cdot 23 \cdot 61$ & $-2^{18}
\cdot 166712087$\\
& $3$ & $-2^3 \cdot 3^6 \cdot 2328401$ & $-3^{12} \cdot 19 \cdot 6547
\cdot 40968624383$\\
&$5$&$2^2 \cdot 3^2 \cdot 5^3 \cdot 1562781531383$&$5^6 \cdot 29 \cdot
37 \cdot 7793 \cdot 31534787 $\\& & &$\cdot
2826173488483$\\
\hline
$\Upsilon 26a$&$2$&$-2^{13} \cdot 3^2 \cdot 5 \cdot 7^2$&$-2^{26} \cdot
7 \cdot 19 \cdot 31 \cdot 1493$\\
&$3$&$-2^3 \cdot 3^5 \cdot 5 \cdot 307 \cdot 61091$&$-3^{10} \cdot
7102940247697920959$\\
&$5$&$-2^2 \cdot 3^2 \cdot 5^5 \cdot 13 \cdot 37 \cdot 293 $&$-5^{10} \cdot 31 \cdot 83 \cdot 1229$\\
& & $\cdot 1847 \cdot 3067$&$ \cdot 155818729039703554943$\\ 
\hline
$\Upsilon 26b$&$2$&$-2^9 \cdot 3^2 \cdot 5 \cdot 229$&$2^{18} \cdot
508045441$\\

&$3$&$-2^3 \cdot 3^7 \cdot 5 \cdot 7 \cdot 1061 \cdot
1579$&$3^{14} \cdot 259103249 \cdot 297496289$\\

&$5$&$2^2 \cdot 3^2 \cdot 5^3 \cdot 7 \cdot 37 \cdot 757$&$-5^6 \cdot 70949\cdot 2914444459 $\\
& &$\cdot 2713
\cdot 51713$&$\cdot
3212880792034321$\\
\hline 
\end{tabular}
\caption{Table of Eigenvalues for genus 2, \cite{Skoruppa}}\label{evalues}
\end{table}

Using Table \ref{evalues} we can compute the Satake $p$-parameters for
the $\Upsilon*$ and primes $p$ in the table. 

First, from the eigen value $T(p^2)$ we need to recover the eigenvalue for $T_i(p^2)$ where $0\leq i\leq 2$.  By an elementary
computation we know 
\[
\lambda_F(T_2(p^2))=p^{2k-6}.
\]
We also know that
\begin{align}\label{eqn1-for-Sat-par} 
\lambda_F(T_0(p))^2&=\left(
  \alpha_0(1+\alpha_1)(1+\alpha_2)\right)^2\nonumber\\
&=\alpha_0^2((1+\alpha_1^2+\alpha_2^2+\alpha_1^2\alpha_2^2)+\nonumber\\
&2(\alpha_1+\alpha_2+\alpha_1^2\alpha_2+\alpha_1\alpha_2^2)+4\alpha_1\alpha_2)\\
&= p^{2k-3}c_1+2c_2p^{2k-3}+4p^{2k-3}\nonumber
\end{align}
\noindent using the notation of system \ref{sys4}.

For $n=2$ we compute the constants in \propref{Krieg} and with them
find the image of $T(p^2)$
\begin{align}\label{eqn2-for-Sat-par}
\lambda_F(T(p^2))&=[\Omega]_2\begin{pmatrix}1\\1\\1\end{pmatrix}\nonumber\\
&=\alpha_0^2((\frac{1}{p^3}+\frac{p^2-1}{p^3}+\frac{2p-2}{p})\alpha_1\alpha_2+\nonumber\\
&(\frac{1}{p}+\frac{p-1}{p})(\alpha_1+\alpha_2+\alpha_1^2\alpha_2+\alpha_1\alpha_2^3)+(1+\alpha_1^2+\alpha_2^2+\alpha_1^2\alpha_2^2))\\
&=p^{2k-4}+(p^2-1)c_2p^{2k-6}+(2p-2)p^{2k-4}+(p-1)c_2p^{2k-4}+c_1p^{2k-3}\nonumber.
\end{align}

Using Maple 9 we solved equations \ref{eqn1-for-Sat-par} and
\ref{eqn2-for-Sat-par} for the constants $c_1$ and $c_2$ to plug into
system \ref{sys4} for $n=2$:
\begin{align*}
f_1:=\sum_{k=0}^2\mathfrak{S}^{(2)}_{k,0}-c_1x_1x_2=0\\
f_2:=\sum_{k=0}^2\mathfrak{S}^{(2)}_{k,1}-c_2x_1x_2=0.
\end{align*}  
Then, we applied the algorithm with the constants
$c_1$ and $c_2$ we just computed and were able to determine Table
\ref{tab-Sat-par} of Satake $p$-parameters.

\begin{table}
\begin{tabular}{|l|r|r|r|}
\hline
$\Upsilon *$ & $p$ & $\alpha_{1,p}$ & $\alpha_{2,p}$\\
\hline
$\Upsilon 20$ & $2$ & $0.6480\dots+i0.7616\dots$ & $-0.2194\dots+i0.9756\dots$\\
& $3$ & $0.4600\dots+i0.8879\dots$ & $-0.9542\dots+i0.2990\dots$\\
& $5$ & $0.6889\dots+i0.7248\dots$ & $-0.9443\dots+i0.3290\dots$\\
& $7$ & $0.1871\dots+i0.9823\dots$ & $-.0922\dots+i0.3865\dotsi$\\
\hline
$\Upsilon 22$ & $2$ & $0.2346\dots+i0.9720\dots$& $-0.5479\dots+i0.8364\dots$\\
& $3$ & $-0.1459\dots+i0.9892\dots$ & $-0.9281\dots+i0.3721\dots$\\
& $5$ & $-0.0257\dots+i0.9996\dots$ & $-0.9532\dots+i0.3022\dots$\\
\hline
$\Upsilon 24a$&$2$&$0.0757\dots+i0.9971\dots$& $-0.7957\dots+i0.6055\dots$\\
&$3$&$0.1368\dots+i0.9905\dots $&$-0.7907\dots+i0.6121\dots$\\
& $5$ &$-0.1066\dots+i0.9943\dots $&$ -0.9999\dots+i0.0046\dots$\\
\hline
$\Upsilon 24b$&$2$& $0.8652\dots+i0.5013\dots $ & $-0.8407\dots+i0.5413\dots$\\
& $3$ & $0.3534\dots+i0.9354\dots $ & $-0.9884\dots+i0.1514\dots $\\
&$5$&$0.2143\dots+i0.9767\dots $&$-0.6417\dots+i0.7668\dots $\\
\hline
$\Upsilon 26a$&$2$&$0.4266\dots+i0.9044\dots $&$-0.8984\dots+i0.4391\dots $\\
&$3$&$0.7854\dots+i0.6189\dots $& $-0.9805\dots+i0.1962\dots $\\
&$5$&$0.1757\dots+i0.9844\dots $&$-0.9034\dots+i0.4287\dots $\\
\hline
$\Upsilon 26b$&$2$&$-0.1533\dots+i0.9881\dots$&$-0.1533\dots+i0.9881\dots $\\

&$3$&$0.5703\dots+i0.8213\dots $&$-0.2998\dots+i0.9539\dots$\\

&$5$&$0.4852\dots+i0.8743\dots $&$-0.8548\dots+i0.5189\dots $\\

\hline 
\end{tabular}
\caption{Table of Satake $p$-parameters for genus 2}\label{tab-Sat-par}
\end{table}

\subsection{Satake $p$-parameters for genus $4$}
The Ikeda lift \cite{Ikeda} is a generalization of the Saito-Kurokawa
lift in that it takes a classical Hecke eigenform $f$ of weight $2k$
to a Hecke eigenform $F$ of weight $k+n$ of genus $n$ (here $k\cong
n\pmod{2}$ whose standard $L$-function is equal to
\[
\zeta(s)\prod_{i=1}^{2n}L(s+k+n-i,f).
\]
It is expected that Hecke eigenforms not in the image of the Ikeda
lift should satisfy the Ramanujan-Petersson conjecture, i.e., should
have unimodular Satake $p$-parameters.  The Schottky form
\cite{Breul2} $J$ is the lifted image of the unique normalized
classical cusp form of weight 12 and is the unique normalized cusp
form of weight $8$ and genus $4$.  Breulmann and Kuss were able to
compute Hecke eigenvalues for $J$; their data is contained in Table
\ref{evalues-J}.  In Table \ref{tab-Sat-par-J} we provide the Satake
$p$-parameters of $J$ based on their computations; these data verify
the expectation described above.

\section{An Application to the Ramanujan-Petersson Conjecture}

As noted above, the Ramanujan-Petersson conjecture for genus 2 claims that a
Siegel modular form $F$ not in the Maass space has Satake
parameters of modulus 1.  This claim has been verified in
Table \ref{tab-Sat-par}.

As noted in \corref{Pn-palindromic}, the output of algorithm
$\mathcal{A}$ is a palindromic polynomial of degree $2n$.  Palindromic
polynomials have reciprocal roots and often have roots on the unit
circle as the following proposition shows:
\begin{prop}[\cite{KonvalinaMatache}]\label{palindromic-roots}  Let $n$ be an even integer and let $p(x)=\sum_{k=0}^{n}
  a_kx^k$ be a palindromic polynomial.  If there exists a
  $k\in\{0,1,\dots,n/2-1\}$ such that 
\[
|a_k|\geq |a_{n/2}|\cos\left (\frac{\pi}{\lfloor
 \frac{n/2}{n/2-k}\rfloor+2}\right)
\]
then $f(x)$ has a pair of unimodular roots.
\end{prop}

In terms of the constants $c_1$ and $c_2$ defined in system
\ref{sys4}, our polynomial $P_2$ is 
\begin{equation}\label{p2}
P_2=x_2^4-c_2x_2^3+(2+c_1)x_2^2-c_2x_2+1.
\end{equation}
\noindent We have the following proposition:
\begin{prop}\label{laurent}  Suppose that $F(x)\in\bbC[x]$ is a palindromic
  polynomial of degree $2n$.  Then
\[
\frac{1}{x^n}F(x)=\overline{F}(x+x^{-1})
\]
where $F(x)\in\bbC[x+x^{-1}]$ is a polynomial of degree $n$.
\end{prop}
\begin{proof}  Let $F(x)=\sum_{i=0}^{2n} a_i x^i$ so that 
\begin{align*}
F(x)/x^n&=\sum_{i=-n}^{n} a_i x^i\\
&= \sum_{i=0}^n a_i (x^i+x^{-i})^.
\end{align*}
We claim that each $x^i+x^{-i}$ is a polynomial
of degree $i$ in $\bbC[x+x^{-1}]$.
\begin{proof}[Proof of Claim]
We prove the statement by induction.  When $i=1$:  clearly
$x+x^{-1}\in\bbC[x+x^{-1}]$.  Assume the statement holds for $i$ and
deduce it holds for $i+1$:
\[
(x^{i+1}+x^{-(i+1)})=(x^i+x^{-i})(x^1+x^{-1})-(x^{i-1}+x^{-(i-1)}).
\]
By induction hypothesis, $(x^i+x^{-i})\in \bbC[x+x^{-1}]$ is of degree
$i$ and $(x^{i-1}+x^{-(i-1)})\in \bbC[x+x^{-1}]$ is of degree $i-1$.
Thus the right hand side is a polynomial of degree $i+1$ in $\bbC[x+x^{-1}]$. 
\end{proof}
In light of the claim, the polynomial $F(x)/x^n\in\bbC[x+x^{-1}]$ is
of degree $n$.
\end{proof}

\begin{table}
\begin{tabular}{|l|r|r|r|r|}
\hline
$p$&$\lambda_J(T_0(p))$&$\lambda_J(T_1(p^2))$&$\lambda_J(T_2(p^2))$&$\lambda_J(T_3(p^2))$\\\hline
$2$&$2^6\cdot 3^3\cdot 5$&$-2^{13}\cdot 3\cdot 5$&$2^{14}\cdot 3^2\cdot 5\cdot
7$&$-2^{14}\cdot 3^3\cdot 5^2$\\\hline
$3$& $2^6\cdot 3^4\cdot 5 \cdot$&$2^4\cdot 3^{11}\cdot $&$2^4\cdot
3^{10}\cdot 5\cdot$&$2^7\cdot 3^9\cdot 5^2\cdot$\\
& $7\cdot 17$&$ 5\cdot 17$&$7\cdot 13\cdot 107$&$7\cdot 1979$\\\hline
$5$&$2^4\cdot 3^3 \cdot 5^2\cdot$&$2^3\cdot 3\cdot
5^{10}\cdot$&$2^3\cdot 3^2\cdot 5^8\cdot$&$2^5\cdot 3^3\cdot 5^6\cdot $\\
   &$7\cdot 131\cdot 199$&$13^2\cdot 41$&$7\cdot 13\cdot 31\cdot
   37253$&$7\cdot 13\cdot 164496949 $\\\hline
$7$&$2^8\cdot 5^2\cdot 7^2\cdot $&$-2^5\cdot 5^2\cdot 7^{10}\cdot $&$2^5\cdot 3 \cdot 5^2\cdot $&$-2^9\cdot 5^4\cdot 7^6\cdot $\\
&$17^2\cdot 1051 $&$ 1049 $&$7^8\cdot 19\cdot659\cdot
23039 $&$1446422309 $\\\hline
\end{tabular}
\caption{Table of Hecke eigenvalues for the Schottky form, \cite{Breul2}}\label{evalues-J}
\end{table}

\begin{table}
\begin{tabular}{|l|r|r|}
\hline
$p$&$\alpha_{1,p}\,,\alpha_{2,p}$&$\alpha_{3,p}\,,\alpha_{4,p}$\\\hline
$2$&$-0.1875\dots\pm i0.6817\dots$ &$-0.7500\dots\pm
i2.7271\dots$\\\hline
$3$&$0.5185\dots\pm i1.6526\dots$&$0.0576\dots\pm
i0.1836\dots$\\\hline
$5$&$0.1545\dots\pm i0.4196\dots$&$0.0309\dots\pm i0.0839\dots$\\\hline
$7$&$-0.4981\dots\pm i2.5984\dots$&$-0.0101\dots\pm i0.0530\dots$\\\hline
\end{tabular}
\caption{Table of Satake $p$-parameters of $J$}\label{tab-Sat-par-J}
\end{table}

The polynomial $P_2$ in equation \ref{p2} can be rewritten as in
\propref{laurent}: 
\[
\frac{1}{x_2^2}P_2(x_2)=\overline{P}_2(x_2+x_2^{-1})=(x_2+x_2^{-1})^2-c_2(x_2+x_2^{-1})+c_1.
\]
Let $\gamma_1,\gamma_2$ be the two roots of
$\overline{P}_2(x_2+x_2^{-1})$.  Note that $P_2$ has four unimodular
roots when both
\[
x_2+\frac{1}{x_2}=\gamma_1 \text{ and } x_2+\frac{1}{x_2}=\gamma_2
\]
have a pair of unimodular roots.  By \propref{palindromic-roots} this
happens when $|\gamma_i|<2$. 

\section{Appendix}  In this appendix we show how to compute explicitly
$\lambda_F(T_0(p^2))$ knowing $\lambda_F(T_0(p))$ and
$\lambda_F(T_i(p^2))$ (for $1\leq i\leq n$).  From this we immediately
see that the implementation of the algorithm requires knowledge of the
eigenvalues of $F$ with respect to the generators and that knowledge
of $\lambda_F(T_0(p^2))$ is redundant.

First, recall the last equations of system \ref{sys2}:
\[
\lambda_{F,0}(p^2)=\sum_{i=1}^{n}
c_{i,n+1}k_i+c_{n+1,n+1}\Sym\left(x_0^2x^{(2,\dots,2)}\right)
\]
\noindent where the constants $c_{i,n+1}$ are known by \propref{Krieg}
and the constants $k_i$ are determined by knowing that
$k_1=p^{nk-n^2-n}$ and the recursion relation:
\[
k_i:=\lambda_F(T_{n-i+1}(p^2)-\sum_{j=1}^{i-1} c_{j,i+1}k_j.
\]
\noindent  So if we can express
$\Sym\left(x_0^2x^{(2,\dots,2)}\right)$ in terms of the constants
$k_i$ we will have shown what we intend to show.  The crux is the
following proposition:
\begin{prop}
\begin{align}\label{propeqn}
\left
  (\lambda_F(T_0(p))\right)^2&=\left(x_0(1+x_1)\cdots(1+x_n)\right)^2\\
&=\sum_{i=0}^n 2^i \Sym_{W_n}\left((x_0^2 x^{\b_i}\right)\nonumber
\end{align}
\noindent where the $\b_i$ we defined at the beginning of \secref{Kriegsec}.
\end{prop}
\begin{proof}  The proof is by induction on $n$.  When $n=1$, the
  left-hand side of system \ref{propeqn} is
\begin{align*}
\sum_{i=0}^1 2^i \Sym_{W_1}\left((x_0^2
  x^{\b_i}\right)&=\Sym_{W_1}\left((x_0^2
  x^{\b_0}\right)+2\Sym_{W_1}\left((x_0^2 x^{\b_1}\right)\\
&=x_0^2(x_1^2+1)+2x_0^2x_1\\
&=(x_0(1+x_1))^2.
\end{align*}
Assume the left-hand side of the equation holds for $n$ and show that
is holds for $n+1$.  Now,
\[
\sum_{i=0}^{n+1} 2^i \Sym_{W_{n+1}}\left((x_0^2 x_1^2\cdots
  x_{n-i+1}^2x_{n-i+2}\cdots x_{n+1}\right)=\sum_{i=0}^{n+1} 2^i\left(\sum_{k=0}^{n+1-i}\mathfrak{S}_{k,i}^{(n+1)}\right)
\]
\noindent by \lemref{frakS}.  Similarly to how we derived equation \ref{f-i-1}
we get the following:
\begin{align*}
\sum_{i=0}^{n+1} 2^i\left(\sum_{k=0}^{n+1-i}\mathfrak{S}_{k,i}^{(n+1)}\right)
&=\sum_{i=0}^{n+1} 2^i\sum_{k=0}^{n+1-i}
x_{n+1}^2\mathfrak{S}_{k-1,i}^{(n)}+x_n\mathfrak{S}_{k,i-1}+\mathfrak{S}_{k,i}\\
&=\sum_{i=0}^{n+1}2^i\sum_{k=0}^{n-i}\mathfrak{S}_{k,i}(x_{n+1}^2+1)+\sum_{i=1}^{n+1}2^ix_{n+1}\sum_{k=0}^{n+1-i}\mathfrak{S}_{k,i-1}^{(n)}\\
&=(x_{n+1}^2+1)\sum_{i=0}^n2^i[x_1^2\cdots x_{n-i}^2x_{n-i+1}\cdots x_n]+\\
&\sum_{i=0}^n
2^{i+1}x_{n+1}\sum_{k=0}^{n-i}\mathfrak{S}_{k,i}^{(n)}\\
&=(x_{n+1}^2+1)\sum_{i=0}^n2^i[x_1^2\cdots
x_{n-i}^2x_{n-i+1}\cdots x_n]+\\&
2x_{n+1}\sum_{i=0}^n2^i[x_1^2\cdots
x_{n-i}^2x_{n-i+1}\cdots x_n]\\
&=(x_{n+1}+1)^2\left ( \sum_{i=0}^n2^i[x_1^2\cdots
x_{n-i}^2x_{n-i+1}\cdots x_n]\right)\\
&= (x_0^2(1+x_1)\cdots(1+x_{n+1}))^2.
\end{align*}
\end{proof}
Rephrasing this in terms of known constants, we have
\begin{align*}
(\lambda_F(T_0(p)))^2&=\sum_{i=0}^n [x_0^2x_1^2\cdots
x_{n-i}^2x_{n-i+1}\cdots x_n]\\
&=[x_0^2x_1^2\cdots x_n^2]+\sum_{i=1}^n
2^i[x_0^2x_1^2\cdots x_{n-i}^2x_{n-i+1}\cdots x_n]\\
&=[x_0^2x_1^2\cdots x_n^2]+\sum_{i=1}^n 2^i k_i.
\end{align*}

So, on the one hand 
\[
[x_0^2x_1^2\cdots x_n^2]=\lambda_F(T_0(p))^2-\sum_{i=1}^n
2^i k_i
\]
and on the other hand by the last equation of system \ref{sys4}
\[
[x_0^2x_1^2\cdots
x_n)^2]=\frac{\lambda_F(T_0(p^2))-\sum_{i=1}^{n}
  c_{i,n+1}k_i}{c_{n+1,n+1}}.
\]
From this it is clear that we can compute $\lambda_F(T_0(p^2))$, and
thus, it suffices, for the algorithm, to know the eigenvalues of $F$
with respect to the generators of the local Hecke algebra.

\bibliographystyle{amsplain}
\bibliography{ryan}

\end{document}